\documentclass[11pt]{article}
\usepackage{amssymb, latexsym, mathtools, amsthm}
\begingroup
    \makeatletter
    \@for\theoremstyle:=definition,remark,plain\do{%
        \expandafter\g@addto@macro\csname th@\theoremstyle\endcsname{%
            \addtolength\thm@preskip\parskip
            }%
        }
\endgroup
\usepackage{enumerate}
\usepackage{pgf,tikz}
\usepackage{pdfsync}
\usepackage{txfonts}
\usepackage[T1]{fontenc}
\usepackage{booktabs}
 \usepackage{graphicx}
\definecolor{dnrbl}{rgb}{0,0,0.3}
\definecolor{dnrgr}{rgb}{0,0.3,0}
\definecolor{dnrre}{rgb}{0.5,0,0}
\usepackage[colorlinks=true, citecolor=dnrgr, linkcolor=dnrre, urlcolor=dnrre] {hyperref}
\usepackage{xcolor}
\usepackage{parskip}
\usepackage{tabularx, colortbl}
\usepackage{geometry}
 \usepackage[scriptsize, up]{caption}
\usepackage{pgf,tikz}
\usetikzlibrary{decorations, decorations.pathmorphing}
\usetikzlibrary {shapes}

\theoremstyle{plain}
\newtheorem{thm}{Theorem}[section]
\newtheorem{prop}[thm]{Proposition}
\newtheorem{lem}[thm]{Lemma}
\newtheorem{coro}[thm]{Corollary}

\numberwithin{equation}{subsection}
\makeatletter

\let\c@table\c@figure
\makeatother

\newcommand{\Nat}{\mathbb{N}}

\newcommand{\LLc}{\mathcal{L}}
\newcommand{\RRc}{\mathcal{R}}

\newcommand{\FSW}{Figueira, Stephan, and Wu\ }

\newcommand{\KS}{Ku{\v{c}}era and Slaman\ }
\newcommand{\CHKW}{Calude, Hertling, Khoussainov and Wang\ }
\newcommand{\DHN}{Downey, Hirschfeldt and Nies\ }
\renewcommand{\DH}{Downey and Hirschfeldt\ }

\newcommand{\ml}{Martin-L\"{o}f }

\newcommand{\ie}{i.e.\ }
\newcommand{\ce}{c.e.\ }
\newcommand{\dce}{d.c.e.\ }
\newcommand{\lce}{left-c.e.\ }
\newcommand{\rce}{right-c.e.\ }

\newcommand{\pf}{prefix-free }

\renewenvironment{abstract}
 { \normalsize
  \list{}{
    \setlength{\leftmargin}{.0cm}%
    \setlength{\rightmargin}{\leftmargin}%
    }%
  \item {\bf \abstractname.} \relax}
 {\endlist}

\title{A note on the differences of computably enumerable reals
\thanks{Barmpalias was supported by the 
1000 Young Talents Plan from the Chinese Government, grant no.\ D1101130.
Additional support was received by
the Chinese Academy of Sciences (CAS) and the Institute of Software of the CAS.
Lewis-Pye was supported by a Royal Society University 
Research Fellowship.}}

\author{George Barmpalias  \and Andrew Lewis-Pye}
\date{\today}
\begin{document}
\maketitle
\begin{abstract}
We show that given any non-computable \lce real $\alpha$ there exists a \lce real $\beta$
such that $\alpha\neq\beta+\gamma$ for all \lce reals and all \rce reals $\gamma$. The proof is 
non-uniform, the dichotomy being whether the given real $\alpha$ is \ml random or not.
It follows that given any universal machine $U$, there is another universal machine $V$ such that
the halting probability $\Omega_U$ of $U$ is not a translation of the halting probability $\Omega_V$ 
of $V$ by a \lce real. We do not know if there is a uniform proof of this fact.
\end{abstract}
\vspace*{\fill}
\noindent{\bf George Barmpalias}\\[0.5em]
\noindent
State Key Lab of Computer Science, 
Institute of Software, Chinese Academy of Sciences, Beijing, China.
School of Mathematics, Statistics and Operations Research,
Victoria University of Wellington, New Zealand.\\[0.2em] 
\textit{E-mail:} \texttt{\textcolor{dnrgr}{barmpalias@gmail.com}}.
\textit{Web:} \texttt{\textcolor{dnrre}{http://barmpalias.net}}\par
\addvspace{\medskipamount}\medskip\medskip
\noindent{\bf Andrew Lewis-Pye}\\[0.5em]  
\noindent Department of Mathematics,
Columbia House, London School of Economics, 
Houghton Street, London, WC2A 2AE, United Kingdom.\\[0.2em]
\textit{E-mail:} \texttt{\textcolor{dnrgr}{A.Lewis7@lse.ac.uk.}}
\textit{Web:} \texttt{\textcolor{dnrre}{http://aemlewis.co.uk}} 

\vfill \thispagestyle{empty}
\clearpage

\section{Introduction}\label{qEDABw7TRV}
The reals which have a computably enumerable left or right Dedekind cut, also known as \ce reals,
play a ubiquitous role in computable analysis and algorithmic randomness.
The differences of  \ce reals, also known as \dce reals, 
form a field under the usual addition and multiplication, as was
demonstrated by Ambos-Spies, Weihrauch, and Zheng 
\cite{Ambos.ea:00}. Raichev \cite{Raichev:05} and Ng \cite{Ng06} showed that this field is real-closed.
Downey, Wu and Zheng \cite{mlq/DowneyWZ04} studied the Turing degrees of \dce reals.
Clearly \dce reals are $\Delta^0_2$ since they can be computably approximated.
Downey, Wu and Zheng \cite{mlq/DowneyWZ04} showed that
every real which is truth-table reducible to the halting problem is Turing equivalent to
a \dce real. However they also showed that there are $\Delta^0_2$ degrees which
do not contain any \dce reals. In this strong sense, \dce reals form a strict subclass of the 
$\Delta^0_2$ reals.  

Despite this considerable body of work on \dce reals, the following rather basic 
question does not have an answer in the current literature. Given a non-computable \ce real $\alpha$, is there a \ce real $\beta$
such that $\alpha-\beta$ is not a \ce real? The answer is, perhaps unsurprisingly, positive. 
We say that a real is \lce or \rce if its left or right Dedekind cut respectively is computably enumerable.
\begin{thm}\label{GRMMMGordg}
If  $\alpha$ is a non-computable \lce real there exists a \lce real $\beta$
such that $\alpha\neq\beta+\gamma$ for all \lce and all \rce reals $\gamma$.
\end{thm}
An interesting aspect of Theorem \ref{GRMMMGordg} is that its proof depends
crucially on the well-developed theory of \ml random \lce reals, and
in particular the methodology developed by
 \DHN in \cite{Downey02randomness}. The proof is 
nonuniform and one has to consider separately the case where $\alpha$ is
\ml random and the case where it is not. 
We do not know if there is a uniform proof of Theorem \ref{GRMMMGordg}, in the sense that
from a \lce approximation to a non-computable real $\alpha$ we can compute a \lce 
approximation to a real $\beta$ such that 
$\alpha\neq\beta+\gamma$ for all \lce and all \rce reals $\gamma$.

Let us focus on the connection with the theory of \ml random \lce reals, as it is crucial in both of the two cases.
It follows from the work of \DHN \cite{Downey02randomness} that: 
\begin{equation}\label{tdLHGAW8Vq}
\parbox{13cm}{if $\alpha,\beta$ are \lce reals
and $\alpha$ is \ml random while $\beta$ is not, then $\alpha-\beta$ is a \ml random \lce real.}
\end{equation}
This, in particular, means that
in Theorem \ref{GRMMMGordg}, $\alpha$ is \ml random if and only if $\beta$ is \ml random.
Moreover we can use this fact in order to reduce
Theorem \ref{GRMMMGordg} to the following special case, 
which we prove in Section \ref{WBEVWr7D1e}.
\begin{lem}\label{gXJTCy7mts}
If $\alpha$ is a \lce real which is neither computable nor \ml random,
 then there exists a \lce real $\beta$
(also not \ml random) such that $\alpha-\beta$ is neither a \lce real nor a \rce real.
\end{lem}
Let us now see how  Theorem \ref{GRMMMGordg} can be derived from this special case.
First, assume that the given $\alpha$ is \ml random. Lemma \ref{gXJTCy7mts} implies the
existence of two \lce reals $\delta_0,\delta_1$ which are not \ml random and
such that $\delta:=\delta_0-\delta_1$ is neither a \lce nor a  \rce real.
Indeed, we can start with any non-computable \lce real $\delta_0$ which is not \ml random 
(such as the halting problem) and apply
Lemma \ref{gXJTCy7mts} in order to get $\delta_1$ with the required properties. Note that
$\delta_1$ is necessarily not \ml random, because otherwise, given that $\delta_0$ is not \ml
random, it would follow from  \eqref{tdLHGAW8Vq} that  $\delta_0-\delta_1$ would be a \rce real.
To establish Theorem \ref{GRMMMGordg} for this case, we choose $\beta=\alpha+\delta$.
First note that $\alpha-\beta$ is not a \lce real or a \rce real, by the choice of $\delta$.
Second, $\beta=(\alpha-\delta_1)+\delta_0$ and $\alpha-\delta_1$
is \ml random by \eqref{tdLHGAW8Vq}, since $\alpha$ is \ml random. Then 
$\beta$ is a \ml random \lce real as the sum of a \ml random \lce real and another \lce real
(a result that was originally proved by Demuth \cite{Dempseu}).
The case of Theorem \ref{GRMMMGordg} when $\alpha$ is not \ml random is exactly
Lemma \ref{gXJTCy7mts}. We note that, as will become apparent in Section \ref{WBEVWr7D1e},
the proof of this case also makes essential use of  \eqref{tdLHGAW8Vq}.

A subclass of the \lce reals are the characteristic functions of \ce sets (viewed as binary expansions).
These reals were called {\em strongly \lce reals} by  \DHN \cite{Downey02randomness}
and are highly non-random reals. It will be clear from the discussion of Section 
\ref{se:prelim} that in Theorem \ref{GRMMMGordg} we cannot (in general) choose
the real $\beta$ to be strongly \lce as in that case, if the given $\alpha$ is \ml random, then 
$\alpha-\beta$ is a \lce real. However the following can be proved using standard finite injury methods.

\begin{prop}[Properly d.c.e.\ reals]\label{gZLwbpfgWh}
There exist strongly \lce reals $\alpha,\beta$ such that $\alpha-\beta$ is
not a \lce real and is not a \rce real. 
\end{prop}

We conclude this discussion with a corollary of Theorem \ref{GRMMMGordg}
in terms of halting probabilities.
The cumulative work of Solovay \cite{Solovay:75},
\CHKW \cite{Calude.Hertling.ea:01} and \KS \cite{Kucera.Slaman:01}
has  shown that the \ml random \lce reals are exactly the halting probabilities of
universal machines. This class remains the same  whether
we consider \pf machines or plain Turing machines.
Here we consider Turing machines operating on strings, and given
an effective list of all Turing machines 
$(M_e)$, a Turing machine $U$ is called {\em universal} if there exists a computable function 
$e\mapsto\sigma_e$ from numbers
to strings such that $U(\sigma_e\ast \tau)=M_e(\tau)$ for all $e$ and all
strings $\tau$. A similar definition applies to universal \pf machines, restricted to 
Turing machines with \pf domain.

Halting probabilities, or equivalently \ml random \lce reals, are all similar in the sense that they all have the same degree with respect to a wide variety of degree structures 
(see \DH \cite[Chapter 9]{rodenisbook}).
A number of results  have been established, however,  which show that halting probabilities may 
differ in certain ways, depending on the universal machine used. For example, 
\FSW \cite{jc/FigueiraSW06} 
showed that for each universal machine $U$ there exists universal machine $V$
such that $\Omega_{U}$ and $\Omega_{V}$ 
have incomparable truth-table degrees. 
Their proof consists of considering 
$\Omega_{V}=\Omega_{U}+X$ for a creative set $X$ like the halting problem, 
and then using the fact from \cite{Bennett88logicaldepth,Calude.Nies:nd}
that no \ml random real truth-table computes a creative set. 
Recall that the use of an oracle computation of a set $A$ from a set $B$ is
an upper bound (as a function of $n$) on the largest position in the oracle $B$ 
queried in the computation of the first
$n$ bits of $A$. 
Frank Stephan (see \cite[Section 6]{IOPORT.05678491})
showed that for each universal machine $U$ there exists universal machine $V$
such that $\Omega_{U}$ cannot
compute  $\Omega_{V}$ with use $n+c$ for any constant $c$.
Recently Barmpalias and Lewis-Pye have improved the use-bound in this statement to
$n+\log n$, while they also showed that 
$\Omega_{U}$, $\Omega_{V}$ can be computed from each other with use $n+2\log n$,
for any universal machines $U,V$.
Along these lines, we can formulate
 Theorem \ref{GRMMMGordg} as follows.
\begin{coro}\label{coro:applc2}
For each universal by adjunction machine $U_0$ there exists another
universal by adjunction machine $U_1$ such that
for all \lce and all \rce reals $\beta$ we have $\Omega_{U_0}\neq \Omega_{U_1} +\beta$.
\end{coro}
This shows that halting probabilities are not always translations of the halting probability of
a fixed universal machine by a \lce or a \rce real.

\section{Overview of \ml random \lce reals}\label{se:prelim}
Some familiarity with the basic concepts of algorithmic information theory and the
basic methods of computability theory would be helpful for the reader. For such background
we refer to one of the monographs \cite{Li.Vitanyi:93, rodenisbook, Ottobook}, 
where the latter two are more focused on computability theory
aspects of algorithmic randomness.
The theory of \lce reals has grown into a significant part of modern algorithmic randomness, and is
best presented in \cite[Chapters 5 and 9]{rodenisbook}.
The present section is an original presentation of some facts regarding \ml random reals
that stem from \cite{Solovay:75,Calude.Hertling.ea:01,Kucera.Slaman:01} and are further
elaborated on in \cite{Downey02randomness}, which are essential for the proof of Theorem
\ref{GRMMMGordg}. Moreover, some of these facts are not given explicitly in the sources above,
but can be recovered from the proofs.

The systematic study of \ml random \ce reals started with 
Solovay in \cite{Solovay:75}, who showed that Chaitin's halting probability of a \pf machine
(a well known \ml random \lce real) has maximum degree
in a degree structure that measures the hardness of approximating \lce reals by
increasing sequences of rationals. This result was complemented by
the work of \CHKW \cite{Calude.Hertling.ea:01} and \KS \cite{Kucera.Slaman:01}, who showed
that these maximally hard to approximate \lce reals are exactly the halting probabilities of universal machines,
which also coincide with the \ml random \lce reals.
The degree structure introduced in \cite{Solovay:75} is now known as
the {\em Solovay degrees of \lce reals} and was extensively studied in
\cite{Downey02randomness}. 
An increasing computable sequence of rationals $(\alpha_i)$ that converges to a real $\alpha$
is called a {\em \lce approximation to $\alpha$}, denoted  $(\alpha_s)\to\alpha$.
The Solovay reducibility $\beta\leq_S\alpha$ 
between \lce reals $\alpha,\beta$ can be defined equivalently by
any of the following clauses:
\begin{enumerate}[\hspace{1cm}(a)]
\item there exists a rational $q$ such that $q\alpha-\beta$ is \lce
\item there exist a rational $q$ and $(\alpha_s)\to\alpha$, $(\beta_s)\to\beta$ such that 
$\beta-\beta_s< q\cdot (\alpha-\alpha_s)$ for all $s$;
\item there exist a rational $q$ and $(\alpha_s)\to\alpha$, $(\beta_s)\to\beta$ such that 
$\beta_{s+1}-\beta_s< q\cdot (\alpha_{s+1}-\alpha_s)$ for all $s$.
\end{enumerate}
Note that the set of rationals $q$ for which one of the above clauses holds is upward closed - if the clause holds for the rational  $q$ then it also holds for all rationals $q'>q$.  
Although it is not explicitly stated in \cite{Downey02randomness}, it follows from the proofs
that when $\beta \leq_S \alpha$,  the infimums of the rationals $q$ for which the clauses (a), (b) and (c) hold are equal.

\KS \cite{Kucera.Slaman:01} proved that: 
\begin{equation}\label{eq:soloredmre}
\parbox{12cm}{if $(\alpha_s)$, $(\beta_s)$ are \lce approximations to 
$\alpha,\beta$ respectively and if $\alpha$ is \ml random, then
$\liminf_s \big[(\alpha-\alpha_s)/(\beta-\beta_s)\big]>0$.}
\end{equation}
In this sense, \ml random \lce reals can only have slow \lce approximations, compared to
any other \lce real and any \lce approximation to it. 
\DHN \cite{Downey02randomness} showed that any \lce approximation to a 
non-random \lce real is considerably
faster than  every \lce approximation to any \ml random real, in the sense that:
\begin{equation}\label{eq:draninf}
\parbox{12cm}{if $(\alpha_s)$, $(\beta_s)$  are \lce approximations
to $\alpha,\beta$ respectively, $\beta$ is \ml random and $\alpha$ 
is not \ml random, then
$\lim_s \big[(\alpha-\alpha_s)/(\beta-\beta_s)\big]=0$.}
\end{equation}
Demuth \cite{Dempseu}
showed that if $\alpha,\beta$ are \lce reals and at least one of them is 
\ml random, then $\alpha+\beta$ is also \ml random.
\DHN \cite{Downey02randomness} proved that the converse also holds, i.e.: 
\begin{equation}\label{eq:d02ranb}
\parbox{13.5cm}{if $\alpha,\beta$ are \lce reals and $\alpha+\beta$ is \ml random then at least one of
$\alpha,\beta$  is \ml random.}
\end{equation}
We conclude our overview with a proof of \eqref{tdLHGAW8Vq} 
which is essential for the proof of Theorem \ref{GRMMMGordg}, but
which is not stated or proved in  \cite{Downey02randomness} (although
it follows easily from the arguments in that paper). We need the following fact
which was proved in \cite{Downey02randomness} (but stated in a weaker form)
and which is also related to the above discussion regarding clauses (a)-(c). 

\begin{lem}[\DHN \cite{Downey02randomness}]\label{le:doranniisrc}
Suppose that $\alpha,\beta$ have \lce approximations $(\alpha_s), (\beta_s)$
such that  $\forall s\ \big(\alpha-\alpha_s< q\cdot (\beta-\beta_s)\big)$ for some rational $q>0$.
If $p>q$ is another rational, then there exists a \lce approximation $(\gamma_s)$ to $\alpha$
such that $\forall s\ \big(\gamma_{s+1}-\gamma_s< p\cdot (\beta_{s+1}-\beta_s)\big)$.
\end{lem}
Now for \eqref{tdLHGAW8Vq}, assume that
$\alpha$ is \ml random and $\beta$ is not \ml random. By
\eqref{eq:draninf} 
for each \lce approximation $(\alpha_s)$ to $\alpha$ 
there exists a \lce approximation $(\beta_s)$ to $\beta$ such that 
$\beta-\beta_s<2^{-1}\cdot (\alpha-\alpha_s)$ for all $s$. Then by
Lemma \ref{le:doranniisrc} there exists a 
 \lce approximation $(\gamma_s)\to\beta$ such that 
$\gamma_{s+1}-\gamma_s<\alpha_{s+1}-\alpha_s$ for all $s$. 
This means that the approximation $(\alpha_s-\gamma_s)$ to $\alpha-\beta$ is
an increasing \lce approximation. So $\alpha-\beta$ is a \lce real.
It remains to show that $\alpha-\beta$ is \ml random. Since $\beta$ is not \ml random, 
by \eqref{eq:d02ranb} it suffices to show that $(\alpha-\beta)+\beta$ is \ml random.
The latter follows from the hypothesis that $\alpha$ is \ml random.

\section{Proof of Lemma \ref{gXJTCy7mts}}\label{WBEVWr7D1e}
We can use a priority injury construction. 
Let $(\gamma_s^i), (\delta_s^i)$ be an effective list of all increasing and decreasing computable
sequences of rationals in $(0,1)$ respectively. Let $\gamma^i$ be the limit of
$(\gamma_s^i)$ and let $\delta^i$ be the limit of $(\delta_s^i)$.
Given $\alpha$ as in the statement of the lemma, 
it suffices to construct a \lce real $\beta$ such that the following conditions are met:
\[
\LLc_i:\ \alpha-\beta\neq \gamma^i
\hspace{0.7cm}\textrm{and}\hspace{0.7cm}
\RRc_i:\ \alpha-\beta\neq \delta^i.
\]
Given an increasing computable sequence of rationals
$(\alpha_s)$ that coverges to $\alpha$, our construction will define an
increasing sequence of rationals $(\beta_s)$  converging to
$\beta$ such that the above requirements are met.
We list the requirements in order of priority as $\LLc_0, \RRc_0,\LLc_1,\dots$.

\paragraph{\bf Parameters of the construction.}
 Let $\beta_0=0$.
The strategy for $\LLc_i$ will use a dynamically defined parameter $c_i$ 
and the strategy for  $\RRc_i$ will use a similar parameter $d_i$. Let $c_i[0]=d_i[0]=0$.
We say that stage $s+1$ is $\LLc_i$-expansionary if 
$|\alpha_{s+1}-\beta_{s+1}-\gamma^i_{s+1}|< 2^{-c_i[s]}$. Similarly, stage $s+1$ is 
$\RRc_i$-expansionary if 
$|\alpha_{s+1}-\beta_{s+1}-\delta^i_{s+1}|< 2^{-d_i[s]}$.
The strategy for each  requirement $\LLc_i$ will define a \lce real
$\beta^i$, which will be its contribution toward the global \lce real $\beta$.
Formally, given the approximations $(\beta^i_s)$ defined by the requirements $\LLc_i$
respectively, for each $s$ we define:
\[
\beta_s=\sum_{i\leq s} \beta^i_s.
\]
If $s+1$ is  $\LLc_i$-expansionary we let $c_i[s+1]=c_i[s]+1$, and otherwise we let
$c_i[s+1]=c_i[s]$. Similarly, if 
$s+1$ is  $\RRc_i$-expansionary we let $d_i[s+1]=d_i[s]+1$, and if not we let
$d_i[s+1]=d_i[s]$. This completes the definition of the parameters $c_i,d_i$ throughout the
stages of the construction.
At each stage $s+1$ the strategy for $\RRc_i$ imposes an automatic restraint on the strategies for $\LLc_j$ of lower priority, which prohibits any increase of $\beta$ by more than $2^{-d_i[s+1]}$.
All of the strategies for the  $\LLc_i$ requirements will use a fixed \ml random \lce real $\eta\in (0,1)$ and an increasing computable
rational approximation $(\eta_s)$ to $\eta$. The  strategy for  each $\LLc_i$ has an extra parameter
$q_i$, which is updated during the stages $s$ and which dictates the scale at which $\eta$ is
going to affect the growth of $(\beta_s)$. At stage $s+1$ we define $q_0[s+1]=\frac{1}{2}$,  and for $i>0$ we define $q_i[s+1]$ to be
the least of all $2^{-i-d_j[s+1]-1}$ for $j<i$.  

\paragraph{\bf Construction of $(\beta_s)$.}  At each stage $s+1$ and each $i\leq s$,
if $s+1$ is  $\LLc_i$-expansionary we define 
$\beta^i_{s+1}=\beta^i_s+ q_i[s+1]\cdot (\eta_{s+1}-\eta_t)$, where
$t$ is the largest $\LLc_i$-expansionary stage before $s+1$ if there is such, and where  $t=0$ otherwise. 
If $s+1$ is not $\LLc_i$-expansionary, we define 
$\beta^i_{s+1}=\beta^i_s$. This completes the definition of $(\beta_s)$.
\paragraph{\bf Verification.}
First we verify that $(\beta_s)$ reaches a finite limit $\beta$. 
Let $\beta^i$ be the limit of $\beta^i_s$ as $s\to\infty$ and note that for each $i$:
\[
\beta^i\leq 2^{-i-1}\cdot \eta < 2^{-i-1}
\hspace{0.5cm}\textrm{so}\hspace{0.5cm}
\beta=\sum_i\beta^i<1.
\]
Recall the dynamic definition of  $c_i[s]$ and $d_i[s]$. It follows that if
$c_i[s]$ reaches a limit, requirement $\LLc_i$ is met. Similarly, if 
$d_i[s]$ reaches a limit, requirement $\RRc_i$ is met. We prove both of these statements by induction.
Suppose that the claim holds for all $i<n$.  Also let $s_0$ be a stage 
such that $c_i[s]=c_i[s_0]$ and $d_i[s]=d_i[s_0]$ for all $i<n$ and all $s>s_0$.
Then by definition $q_n[s]= q_n[s_0]$ for all $s>s_0$. Let $q_n$ denote the limit $q_n[s_0]$ of
$q_n[s]$ from now on.
If $c_n[s]$ does  not reach a limit, then there are infinitely many 
$\LLc_n$-expansionary stages, which implies that $\alpha-\beta=\gamma^n$.
Moreover if $(t_j)$ is a monotone enumeration of the $\LLc_n$-expansionary stages, then
$\beta_{t_{s+1}}-\beta_{t_{s}}>q_n\cdot (\eta_{t_{s+1}}-\eta_{t_{s}})$ for all $s$.
Since $\eta$ is \ml random, this means that $\beta$ is also \ml random. But by hypothesis $\alpha$
is not \ml random, so $\alpha-\beta$ is a \ml random \rce real. This contradicts the fact that
$\alpha-\beta=\gamma$ since \rce reals which have a \lce approximation are computable.
It follows that there are only finitely many $\LLc_n$-expansionary stages, which implies that
$c_n[s]$ reaches a limit. Let $s_1>s_0$ be a stage such that $c_n[s]=c_n[s_1]$ for all $s>s_1$.

It remains to show that $d_n[s]$ reaches a limit. Towards a contradiction, suppose that this is not the case,
so that there are infinitely many $\RRc_n$-expansionary stages.
Then it follows that $\alpha-\beta=\delta^n$. Let $(t_k)$ be a computable enumeration of all
$\RRc_n$-expansionary stages. Then $d_n[t_k]=k$ for all $k$. For each $i>n$ and 
each $k$ we have
$\beta^i-\beta^i_{t_k}\leq 2^{-i-k-1}$ which means that for $k$ large enough that $t_k>s_1$: 
\[
\beta-\beta_{t_k} < \sum_{i>n} (\beta^i-\beta^i_{t_k}) \leq \sum_{i>n} 2^{-i-k-1} \leq  2^{-k-1}.
\]
This means that $\beta$ is a computable real. Since $\alpha=\delta^n+\beta$ and 
$\delta^n$ is a \rce real, it follows that $\alpha$ is a \rce real. Since $\alpha$  also a \lce real, it must therefore be
computable, contrary to hypothesis.  So we may conclude that
there are finitely many $\RRc_n$-expansionary stages, which establishes that $\RRc_n$ is met and
$d_n$ reaches a limit. This concludes the induction step and the proof that the constructed real $\beta$
meets the requirements $\LLc_n$ and $\RRc_n$ for all $n$. 

\paragraph{\bf Remark.}
The reader may wonder why a uniform argument for Theorem \ref{GRMMMGordg} might not work, \ie
why we needed to divide into two cases according to whether the given real is \ml random or not.
While it is not easy to explain why some things do not work, the immediate answer is that in a construction
such as the argument above, if we did not assume that the given real is not \ml random or we did not code
randomness into the real we construct, we do not see a way to argue that requirements $\LLc_i$ act only 
finitely often.
More generally, if a direct standard uniform construction worked, in our view we could use it to show that 
{\em given a left ce real $\alpha$ we can find a left ce real $\beta$ such that $2\alpha-\beta$ is not \lce and 
$\alpha-\beta$ not a \rce real}.
However we know that this is not possible by one of the results in \cite{omegax}. 
This non-uniformity seems to relate to the non-uniformities
in the characterization of the halting probabilities in  \cite{Solovay:75,Calude.Hertling.ea:01,Kucera.Slaman:01}
that we discussed in Section \ref{qEDABw7TRV}.
Showing that such non-uniformities are necessary may be an interesting exercise.

\section{Proof of Proposition \ref{gZLwbpfgWh}}
We can use a standard priority injury construction. 
Let $(\gamma_s^i), (\delta_s^i)$ be an effective list of all increasing and decreasing computable
sequences of rationals in $(0,1)$ respectively. Moreover let $\gamma^i$ be the limit of
$(\gamma_s^i)$ and let $\delta^i$ be the limit of $(\delta_s^i)$.
It suffices to satisfy the following conditions.
\[
\LLc_i:\ \alpha-\beta\neq \gamma^i
\hspace{0.7cm}\textrm{and}\hspace{0.7cm}
\RRc_i:\ \alpha-\beta\neq \delta^i
\]
Our construction will define increasing sequences $(\alpha_s), (\beta_s)$ of rationals which converge
to $\alpha,\beta$ respectively. Let $\alpha_0=\beta_0=0$.
Strategies $\LLc_i$ will use a parameter $c_i$ which takes values from $\Nat^{[2i]}$ (i.e.\ the even numbers)
and strategies $\RRc_i$ will use a parameter $d_i$ which takes values from $\Nat^{[2i+1]}$.
We say that $\LLc_i$ requires attention at stage $s+1$ if either $c_i$ is undefined, 
or $c_i[s]$ is defined and
$|\alpha_s-\beta_s-\gamma^i_{s+1}|< 2^{-c_i[s]-3}$.
Similarly we say that $\RRc_i$ requires attention at stage $s+1$ if
either $d_i$ is undefined, or $d_i[s]$ is defined and
$|\alpha_s-\beta_s-\delta^i_{s+1}|< 2^{-d_i[s]-3}$.
Strategy $\LLc_i$ will impose a restraint $\ell_i$ on $\alpha$ while strategy 
$\RRc_i$ will impose a restraint $r_i$ on $\beta$. The parameters $\ell_i,r_i$ will be defined 
(and possibly redefined) dynamically during the construction, before reaching a limit.
We list the requirements in order of priority as $\LLc_0, \RRc_0,\LLc_1,\dots$  and construct $\alpha,\beta$ 
as \ce sets $A,B$ with characteristic sequences the binary expansions of $\alpha,\beta$. In this way, the
restraints $\ell_i,r_i$ will apply to the enumerations into $A$ and $B$ respectively.
Note that enumerating a number $n$ into $A$ increases $\alpha-\beta$ by $2^{-n}$
while enumerating $n$ into $B$ decreases $\alpha-\beta$ by $2^{-n}$.
Initializing requirement $\LLc_i$ at stage $s+1$ means to let $c_i[s+1], \ell_i[s+1]$ be undefined.
Similarly, initializing  $\RRc_i$ at stage $s+1$ means to let $d_i[s+1], r_i[s+1]$ be undefined.
If $c_i[s]$ is defined and $\LLc_i$ is not initialized at stage $s+1$ then we automatically assume that
 $c_i[s]=c_i[s+1]$. Similarly,
 if $d_i[s]$ is defined and $\RRc_i$ is not initialized 
 at stage $s+1$ then we automatically assume that
 $d_i[s]=d_i[s+1]$.
  
At stage $s+1$ let $i$ be the least number $\leq s$ such that $\LLc_i$ or $\RRc_i$
requires attention. If there is no such number, go to the next stage. Otherwise, first assume  
that $\LLc_i$ requires attention at stage $s+1$. If $c_i[s]$ is not defined, let $c_i[s+1]$ be the least number in
$\Nat^{[2i]}$ which is larger than any value of any parameter defined so far in the construction 
(in particular larger than all previous values of $c_i$ and larger than any restraint $r_j$ on $\beta$ 
which is currently defined). 
If,on the other hand $c_i[s]$ is defined, then enumerate it into $B$, define $\ell_i[s+1]=c_i[s]+3$
and initialize all $\LLc_{j+1}, \RRc_j$ for all $j\geq i$.
In this latter case we say that $\LLc_i$ {\em acts}
at stage $s+1$.

Second, assume  
that $\RRc_i$ requires attention at stage $s+1$. If $d_i[s]$ is not defined, let $d_i[s+1]$ 
be the least number in
$\Nat^{[2i+1]}$ which is larger than any value of any parameter defined so far in the construction 
(in particular larger than all previous values of $d_i$ and larger than any restraint $\ell_j$ on $\alpha$ 
which is currently defined). 
If,on the other hand $d_i[s]$ is defined, then enumerate it into $A$, define $r_i[s+1]=d_i[s]+3$
and initialize all $\LLc_{j}, \RRc_{j}$ for all $j\geq i$.
In this latter case we say that $\RRc_i$ {\em acts}
at stage $s+1$.

The construction defined computable enumerations of the sets $A,B$ which in turn define
computable non-decreasing rational 
approximations $(\alpha_s), (\beta_s)$ to the reals $\alpha,\beta$. Since $A,B$ are \ce and
no \ce set is \ml random, we immediately get that $\alpha,\beta$ are not random. It remains to
show that $\alpha,\beta$ meet the requirements $\LLc_i$ and $\RRc_i$. Note that if 
$\LLc_i$ {\em acts} at stage $s+1$ and is not initialized at any later stage, then it will not 
require attention at any later stage. Indeed, in this case no higher priority requirement will
act at later stages, and both $c_i[t]$ and $\ell_i[t+1]$ remain constant for all $t\geq s$.
Let $c_i,\ell_i$ denote their final values respectively.
Since $\LLc_i$ required attention 
at stage $s+1$ we have $|\alpha_s-\beta_s-\gamma^i_{s+1}|< 2^{-c_i-3}$.
Moreover $\beta_{s+1}-\beta_s=2^{-c_i}$  and $\alpha_s=\alpha_{s+1}$. So
$\alpha_{s+1}-\beta_{s+1}<\gamma^i_{s+1}- 2^{-c_i-1}$
and since $\ell_i=c_i+3$ we have $\alpha_t-\alpha_{s+1}< 2^{-c_i-2}$ for all $t>s$.
Therefore 
$\alpha_{t}-\beta_{t}<\gamma^i_{t}- 2^{-c_i-2}$ for all $t>s$ and  $\LLc_i$ will not require attention
at any stage after $s$.
Moreover we also get that $\alpha-\beta\leq\gamma^i- 2^{-c_i-2}$ which means that
in this case condition $\LLc_i$  is met.
We have shown that:
\begin{equation}\label{zjxyJeyuQK}
\parbox{11cm}{If $\LLc_i$ {\em acts} at stage $s+1$ 
and is not initialized at any later stage, then it will not 
require attention at any later stage and is satisfied.}
\end{equation}
An entirely similar argument shows that: 
\begin{equation}\label{VGMvXj6eIl}
\parbox{11cm}{If $\RRc_i$ {\em acts} at stage $s+1$ 
and is not initialized at any later stage, then it will not 
require attention at any later stage and is satisfied.}
\end{equation}
It remains to use \eqref{zjxyJeyuQK} and \eqref{VGMvXj6eIl} inductively in order to show that
$\alpha-\beta$ meets $\LLc_i$, $\RRc_i$ for all $i$. 
Note that $\LLc_0$ cannot be initialized. So $c_0$ will be defined and remain constant for the
rest of the stages. If $\LLc_0$ never acts, then it does not require attention after the first time that it
required (and received) attention. This means that $|\alpha_s-\beta_s-\gamma^0|\geq 2^{-c_i-3}$ for
all but finitely many stages $s$, so $\alpha-\beta\neq \gamma^i$. If it does act at some stage, then by
\eqref{zjxyJeyuQK}  it is satisfied and never requires attention at any later stage.
Now inductively assume that the same is true for all $\LLc_{i},\RRc_{i}$, $i<e$. Then consider
a stage $s_0$ after which none of $\LLc_{i},\RRc_{i}$, $i<e$ acts or requires attention. Then the same 
argument shows that $\LLc_{e}$ does not act or require attention after a certain stage, and is met.
The same argument applies to $\RRc_{e}$ through property \eqref{VGMvXj6eIl}, 
and this concludes the induction step. We can conclude that $\alpha-\beta$
meets $\LLc_{i},\RRc_{i}$ for all $i$. 

\paragraph{\bf Remark.}
The referee has pointed out that a proof of
Proposition \ref{gZLwbpfgWh} may be given without a direct construction.
Consider two \ce sets $A,B$ such that $A-B$ has properly d.c.e.\ degree, \ie there is no \ce set which is Turing equivalent
to $A-B$. Such \ce sets were originally constructed
in Cooper \cite{Cooper:71}, and the standard construction gives $B\subseteq A$.
Let $\alpha,\beta$ be the reals in $(0,1)$ whose binary expansions are the characteristic sequences of $A,B$ respectively.
Then the binary expansion of $\alpha-\beta$  is the characteristic sequence of $A-B$.
If $\alpha-\beta$ had a \lce or a \rce approximation, then $A-B$ would be Turing equivalent to the left or the right Dedekind cut
of $\alpha-\beta$ which would be a \ce set. This would contradict the choice of $A-B$. Hence $\alpha, \beta$ have
the required properties.


\end{document}